\documentclass[10pt,a4paper]{article}
\usepackage{amsfonts}
\usepackage{amsthm}
\usepackage{amssymb}
\usepackage[centertags]{amsmath}
\usepackage{amssymb}
\usepackage{enumerate}
\usepackage{graphicx}
\providecommand{\U}[2]{\protect\rule{1.5in}{1.5in}}
\theoremstyle{plain}
\newtheorem{theorem}{Theorem}[section]

\newtheorem{definition}[theorem]{Definition}
\newtheorem{example}[theorem]{Example}

\numberwithin{equation}{section}
\begin{document}
\bibliographystyle{plain}
\title{Application of fuzzy Laplace transforms for solving fuzzy partial Volterra integro-differential equations}
\author{Saif Ullah\footnote{Department of Mathematics, University of Peshawar, 25120, Khyber Pakhtunkhwa, Pakistan. E-mail: saifullah.maths@upesh.edu.pk}, Muhammad Farooq\footnote{Department of Mathematics, University of Peshawar, 25120, Khyber Pakhtunkhwa, Pakistan. E-mail: mfarooq@upesh.edu.pk}, Latif Ahmad\footnote{a. Shaheed Benazir Bhutto University, Sheringal. b. Department of Mathematics, University of Peshawar, 25120, Khyber Pakhtunkhwa, Pakistan. E-mail: ahmad49960@yahoo.com}, Saleem Abdullah\footnote{Department of Mathematics, Quaid-i-Azam University, Islamabad, Pakistan. E-mail: saleemabdullah81@yahoo.com}}
%%\begin{document}

\maketitle
\begin{abstract}
Fuzzy partial integro-differential  equations have a major role in the fields of science and engineering. In this paper, we propose  the solution of fuzzy  partial  Volterra integro-differential equation with convolution type kernel using fuzzy Laplace transform  method (FLTM) under  Hukuhara differentiability. It is shown  that FLTM is a simple and reliable approach for solving such equations analytically.  Finally, the method is illustrated with few examples to show the ability of the proposed method.
\end{abstract}
\section{Introduction}
The topic of fuzzy integro differential equations (FIDEs) has been rapidly grown recent years. The   basic idea and arithmetics of fuzzy sets were first introduced by L. A. Zadeh in \cite{21:lt}. The concept of fuzzy derivatives and fuzzy integration were studied in \cite{11:lt, 9:lt} and then some generalization have been investigated in \cite{10:lt, 11:lt, 12:lt, 13:lt}. One of the most important field of the fuzzy theory is the  fuzzy differential equations \cite{1:lt, 2:lt}, fuzzy integral equations \cite{3:lt, 4:lt, 5:lt} and fuzzy integro-differential equations (FIDEs) \cite{6:lt, 7:lt, 8:lt}. The  FIDEs is obtained when a physical system is modeled under differential sense \cite{26:lt}. Also FIDEs in fuzzy setting are a natural way to model uncertainty of dynamical systems. Therefore the solution of the fuzzy integro-differential equations is very important in various fields such as Physics, Geographic, Medical and Biological Sciences \cite{27:lt, 28:lt, 29:lt}. In \cite{9:lt} Seikkala defined fuzzy  derivatives  while concept of integration of fuzzy functions was first introduced by Dubois and Prade \cite{11:lt}. Alternative approaches were later studied in \cite{14:lt, 15:lt}. The idea of fuzzy partial differential equations (FPDEs) was first introduced by Buckley in \cite{31:lt}. Allahveranloo proposed the difference method for solving FPDEs in \cite{32:lt}.

In \cite{18:lt} Allahveranloo and Salahshour  proposed the idea of fuzzy Laplace transform method for solving first order fuzzy differential equations under generalized H-differentiability. The technique of fuzzy Laplace transform method to solve fuzzy convolution Volterra integral equations (FCVIEs) of the second kind  was developed in \cite{19:lt}. Recently the technique used in \cite{19:lt}  was extend  for  solving fuzzy convolution  Volterra integro differential equations (FCVIDEs) in \cite{33:lt} under generalized Hukuhara differentiability.

In \cite{34:lt} the solution of classical  PIDEs was discussed using classical  Laplace transform. In the present  article we investigate the solution of different types of fuzzy partial integro differential equations with convolution  kernel (FPIDEs) using fuzzy Laplace transform method. In order to determine the lower and upper  functions of the solution  we convert the given  FPIDEs to two crisp ordinary differential equations by using FLT.

\noindent Rest of the paper is organized as follows:\\
In section 2, some basic definitions and results are stated which will be used throughout this paper. In section 3, two dimensional fuzzy Laplace transform  is given  and fuzzy convolution theorem is stated in this case. In section 4, the  fuzzy Laplace transform is applied to fuzzy partial Volterra integro-differential equation to construct the general technique. Illustrative examples are also considered to show the ability of the proposed method in section 5, and the conclusion is drawn in section 6.

\section{Preliminaries}
In this section we will recall some basics definitions and theorems needed throughout the paper such as fuzzy number, fuzzy-valued function and the derivative of the fuzzy-valued functions \cite{20:lt, 21:lt}.
%%%\end{preliminaries}
\begin{definition} A fuzzy number is defined  as the mapping such that $u:R\rightarrow[0,1]$, which satisfies the following four properties
\begin{enumerate}
\item $u$ is upper semi-continuous.
 \item $u$ is fuzzy  convex that is $u(\lambda  x+(1-\lambda)y) \geq \min{\{u(x), u(y)\}}.\; x, y\in R$ and $\lambda\in [0,1]$.
\item $u$ is normal that is $\exists$ $x_0\in R$, where $u(x_0)=1$.
\item $A=\{\overline{x \in \mathbb{R}: u(x)>0}\}$ is compact, where $\overline{A}$ is closure of $A$.
\end{enumerate}
\end{definition}
\begin{definition}
A fuzzy number in parametric form is given as an order pair of the form $u=(\underline{u}(r), \overline{u}(r))$, where $0\leq r\leq1$ satisfying the following conditions.
\begin{enumerate}
\item $\underline{u}(r)$ is a bounded left continuous increasing function in the interval $[0,1]$.
\item $\overline{u}(r)$ is a bounded left continuous decreasing function in the interval $[0,1]$.
\item $\underline{u}{(r)\leq\overline{u}(r)}$. % \mbox{ where  }   $0\leq r \leq 1 (though we dont need it as r is already defined above)$
\end{enumerate}
If $\underline{u}(r)=\overline{u}(r)=r$, then $r$ is called crisp number.
\end{definition}

\noindent Since each $y\in R$ can be regarded as a fuzzy number if
\begin{eqnarray*}\widetilde{y}(t)=\begin{cases}1, \;\;\; if \;\; y=t,\\  0, \;\;\; if \;\; y\neq t.\end{cases}\end{eqnarray*}
For arbitrary fuzzy numbers $u=(\underline{u}(\alpha), \overline{u}(\alpha))$ and $v=(\underline{v}(\alpha), \overline{v}(\alpha))$ and an arbitrary crisp number $j$, we define addition and scalar multiplication as:
\begin{enumerate}
\item $(\underline{u+v})(\alpha)=(\underline{u}(\alpha)+\underline{v}(\alpha))$.
\item $(\overline{u+v})(\alpha)=(\overline{u}(\alpha)+\overline{v}(\alpha))$.
\item $(j\underline{u})(\alpha)=j\underline{u}(\alpha)$, $(j\overline{u})(\alpha)=j\overline{u}(\alpha)$, \mbox{       }  $j\geq0$.
\item $(j\underline{u})(\alpha)=j\overline{u}(\alpha)\alpha, (j\overline{u})(\alpha)=j\underline{u}(\alpha)\alpha$, $j<0$.
\end{enumerate}
\begin{definition} (See \cite{18:lt,23:lt}) Let us suppose that x, y $\in E$, if $\exists$ $z\in E$ such that
$x=y+z$, then $z$ is called the H-difference of $x$ and $y$ and is given by $x\ominus y$.\end{definition}

\noindent The Housdorff distance between the fuzzy numbers \cite{6:lt,12:lt,18:lt,23:lt} defined by
\[d:E\times E\longrightarrow R^{+}\cup \{{0}\},\]
\[d(u,v)=\sup_{r\in[0,1]}\max\{|\underline{u}(r)-\underline{v}(r)|, |\overline{u}(r)-\overline{v}(r)|\},\] \noindent where $u=(\underline{u}(r), \overline{u}(r))$ and $v=(\underline{v}(r), \overline{v}(r))\subset R$.% has been utilized by Bede and Gal \cite{12}.
\\\\
We know that if $d$ is a metric in $E$, then it will satisfy the following properties, introduced by Puri and Ralescu \cite{22:lt}:
\begin{enumerate}
\item $d(u+w,v+w)=d(u,v)$, $\forall$  u, v, w $\in$ E.

\item $(k \odot u, k \odot v)=|k|d(u, v)$, $\forall$ k $\in$ R, \mbox{  and  } u, v $\in$ E.

\item $d(u \oplus v, w \oplus e)\leq d(u,w)+d(v,e)$, $\forall$ u, v, w, e $\in$  E.
\end{enumerate}

\begin{theorem} (see Wu \cite{24:lt}) Let $f$ be a fuzzy-valued function on $[a,\infty)$ given in the parametric form as $(\underline{f}(x,r), \overline{f}(x,r))$ for any constant number $r\in[0,1]$. Here we assume that $\underline{f}(x,r)$ and $\overline{f}(x,r)$ are Riemann-Integrable on $[a,b]$ for every $b\geq a$. Also we assume that $\underline{M}(r)$ and $\overline{M}(r)$ are two positive functions, such that
$\int_a^b|\underline{f}(x,r)| dx \leq \underline{M}(r)$ and $\int_a^b |\overline{f}(x,r)| dx \leq \overline{M}(r)$
for every $b\geq a$, then $f(x)$ is improper fuzzy Riemann-integrable on $[{a}, \infty)$. Thus an improper integral will always be a fuzzy number. In short \[ \int_a^b f(x) dx = \bigg( \int_a^b\underline{f}(x,r) dx, \int_a^b \overline{f}(x,r) dx\bigg).\]
%It is will known that Hukuhare differentiability for fuzzy function was introduced by Puri \& Ralescu in \cite{22:lt}. %It is based on H-differentiable.
\end{theorem}
\noindent Next we define the $nth$ order partial H-derivatives for fuzzy valued functions $u=u(x,t)$ with respect to x and t in a similar way as given in \cite{30:lt, 25:lt}.
 \begin{definition} The function  $u:(a,b)\times(a,b)\rightarrow E$, is said to be H-differentiable of the $nth$ order  at   $t_{0}\in (a,b)$, with respect to  $t$, if $\exists$ an element $\frac{\partial^{n} }{\partial t^{n}}u(x,t_0)\in E$ such that
\begin{enumerate}
 \item $\forall h>0$ sufficiently small $\exists$ $\frac{\partial^{n-1}}{\partial t^{n-1}}u(x,t_0+h)\ominus\frac{\partial^{n-1} }{\partial t^{n-1}}u(x,t_0)$, $\frac{\partial ^{n-1}}{\partial t^{n-1}}u(x,t_0)\ominus \frac{\partial^{n-1} }{\partial t^{n-1}}u(x,t_0-h)$, then the following limits hold (in the metric $d$)\\
  $\lim_{h\rightarrow 0}\frac{\frac{\partial^{n-1} }{\partial t^{n-1}}u(x,t_0+h)\ominus\frac{\partial ^{n-1}}{\partial x^{n-1}}u(x,t_0)}{h}=\lim_{h\rightarrow 0}\frac{\frac{\partial ^{n-1}}{\partial t^{n-1}}u(x,t_0)\ominus \frac{\partial^{n-1} }{\partial t^{n-1}}u(x,t_0-h)}{h}=\frac{\partial^{n} }{\partial t^{n}}u(x,t_0)$,
 or
  \item $\forall h>0$ sufficiently small $\exists$ $\frac{\partial ^{n-1}}{\partial t^{n-1}}u(x,t_0)\ominus\frac{\partial^{n-1}}{\partial t^{n-1}}u(x,t_0+h)$, $\frac{\partial ^{n-1}}{\partial t^{n-1}}u(x,t_0-h)\ominus\frac{\partial }{\partial t}u(x,t_0)$, then the following limits hold (in the metric $d$)\\
  $\lim_{h\rightarrow 0}\frac{\frac{\partial^{n-1}}{\partial t^{n-1}}u(x,t_0)\ominus\frac{\partial^{n-1}}{\partial t^{n-1}}u(x,t_0+h)}{-h}=\lim_{h\rightarrow 0}\frac{\frac{\partial ^{n-1}}{\partial t^{n-1}}u(x,t_0-h)\ominus\frac{\partial^{n-1}}{\partial t^{n-1}}u(x,t_0)}{-h}=\frac{\partial^{n} }{\partial t^{n}}u(x,t_0)$.
 \end{enumerate}
 \end{definition}

\noindent Similarly
\begin{definition} The function  $u:(a,b)\times(a,b)\rightarrow E$, is said to be H-differentiable of the $nth$ order  at   $x_{0}\in (a,b)$, w.r.t $x$, if $\exists$ an element $\frac{\partial^{n} }{\partial^{n} x}u(x_0,t)\in E$ such that
\begin{enumerate}
 \item $\forall h>0$ sufficiently small $\exists$ $\frac{\partial^{n-1}}{\partial x^{n-1}}u(x_0+h,t)\ominus\frac{\partial^{n-1} }{\partial x^{n-1}}u(x_0,t)$, $\frac{\partial^{n-1} }{\partial x^{n-1}}u(x_0,t)\ominus \frac{\partial^{n-1} }{\partial x^{n-1}}u(x_0-h,t)$, then the following limits hold (in the metric $d$)\\
  $\lim_{h\rightarrow 0}\frac{\frac{\partial^{n-1} }{\partial x^{n-1}}u(x_0+h,t)\ominus\frac{\partial ^{n-1}}{\partial x^{n-1}}u(x_0,t)}{h}=\lim_{h\rightarrow 0}\frac{\frac{\partial^{n-1} }{\partial x^{n-1}}u(x_0,t)\ominus \frac{\partial ^{n-1}}{\partial x^{n-1}}u(x_0-h,t)}{h}=\frac{\partial^{n} }{\partial x^{n}}u(x_0,t)$,
  or
  \item $\forall h>0$ sufficiently small $\exists$ $\frac{\partial ^{n-1}}{\partial x^{n-1}}u(x_0,t)\ominus\frac{\partial^{n-1}}{\partial x^{n-1}}u(x_0+h,t)$, $\frac{\partial ^{n-1}}{\partial x^{n-1}}u(x_0-h,t)\ominus\frac{\partial ^{n-1}}{\partial x^{n-1}}u(x_0,t)$, then the following limits hold (in the metric $d$)\\
  $\lim_{h\rightarrow 0}\frac{\frac{\partial^{n-1}}{\partial x^{n-1}}u(x_0,t)\ominus\frac{\partial^{n-1}}{\partial x^{n-1}}u(x_0+h,t)}{-h}=\lim_{h\rightarrow 0}\frac{\frac{\partial ^{n-1}}{\partial x^{n-1}}u(x_0-h,t)\ominus\frac{\partial^{n-1}}{\partial x^{n-1}}u(x_0,t)}{-h}=\frac{\partial^{n} }{\partial x^{n}}u(x_0,t)$.

 \end{enumerate}
 The denominators $h$ and $-h$ denote multiplication by $\frac{1}{h}$ and $\frac{-1}{h}$ respectively.
 \end{definition}

 \section{Two dimensional fuzzy Laplace transform}\label{flt}

%\begin{theorem}(See Chalco and Reman-Flores \cite{25:lt}) Let $f:R\rightarrow E$ be a function denoted by

In this section we state some definitions and theorems from \cite{30:lt,19:lt} which will be used in the next section.
\begin{definition}
Let $u=u(x,t)$ is a fuzzy-valued function and $p$ is a real parameter, then  FLT of the function $u$ with respect to t denoted by U(x,p), is defined as follows: \begin{equation*}\label{eq2}U(x,p)=L[u(x,t)]=\int_{0}^{\infty}e^{-pt}u(x,t)dt=\lim_{\tau\rightarrow\infty}\int_{0}^{\tau}e^{-pt}u(x,t)dt,\end{equation*}
\begin{equation*}U(x,p)=\bigg[\lim_{\tau\rightarrow\infty}\int_{0}^{\tau}e^{-pt}\underline{u}(x,t)dt,\lim_{\tau\rightarrow\infty}\int_{0}^{\tau}e^{-pt}\overline{u}(x,t)dt\bigg],\end{equation*}
\noindent whenever the limits exist. The r-cut representation of U(x,p) is given as:
\end{definition}
\begin{equation*}U(x,p;r)=L[u(x,t;r)]=[l(\underline{u}(x,t;r)),l(\overline{u}(x,t;r))],\end{equation*}
\noindent where
\begin{equation*}l[\underline{u}(x,t;r)]=\int_{0}^{\infty}e^{-pt}\underline{u}(x,t;r)dt=\lim_{\tau\rightarrow\infty} \int_{0}^{\tau}e^{-pt}\underline{u}(x,t;r)dt,\end{equation*}
\begin{equation*}\label{eq6}l[\overline{u}(x,t;r)]=\int_{0}^{\infty}e^{-pt}\overline{u}(x,t;r)dt=\lim_{\tau\rightarrow\infty}\int_{0}^{\tau}e^{-pt}\overline{u}(x,t;r)dt.  \end{equation*}
Now to use FLT method we have to state the follwing result:
\begin{theorem}
Let $u:(a,b)\times(a,b)\longrightarrow E$ is a fuzzy valued function  such that its derivatives up to $(n-1)th$ order w.r.t $``t"$ are continuous for all $t>0$  and  $u^{n}$ exists then
\begin{equation}
L[\frac{\partial^n}{\partial x^t}u(x,t)]=p^nU(x,p)\ominus p^{n-1}u(x,0)\ominus p^{n-2}u(x,0)\ominus \cdots \ominus u^{n-1}(x,0)
\end{equation}
\begin{equation}
L[\frac{\partial^n u}{\partial x^n}]=\frac{d^n}{dx^n} L[u(x,t)]=\frac{d^n }{dx^n}U(x,p),
\end{equation}
%and

%\begin{enumerate}
 %\begin{equation}L[\frac{\partial u}{\partial x}]=\frac{\partial}{\partial x}U(x, p),\end{equation}
%\begin{equation}L[\frac{\partial^2 u}{\partial x^2}]=\frac{\partial^2}{\partial x^2}U(x, p),\end{equation}
%and \begin{equation}L[\frac{\partial u}{\partial t}]=pU(x, p)\ominus u(x, 0),\end{equation}
%\begin{equation}L[\frac{\partial^2 u}{\partial t^2}]=p^2U(x, p)\ominus pu(x, 0)\ominus pu'(x, 0),\end{equation}
%\end{enumerate}
\end{theorem}
\begin{definition}
The fuzzy two dimensional  convolution of  fuzzy-valued functions f and g defined by
\begin{equation*}
(f\ast g)(t)=\int_{0}^{t}f(s)g(t-s)ds,
\end{equation*}
where $t>0$ and it exists if f and g are say, piecewise continues functions.
\end{definition}
%\begin{theorem}(Derivative theorem) Suppose that f is continues fuzzy valued function on $[0, \infty)$ and of exponential order \alpha$ and that $f'$ is piecewise continues
% in $[0,\infty)$ then
%\begin{enumerate}
%\item \[L(f'(t))=pL(f(t))\circleddash f(0),\]
%if f is (i)-differentiable.
%\item \[L(f'(t))=(-f(0))\circleddash (-pL(f(t))),\]
%if f is (ii)-differentiable.
%\end{enumerate}
%\end{theorem}
\begin{theorem}\label{ct}(Fuzzy convolution theorem),
if f and  g are piecewise continuous fuzzy-valued function on $[0,\infty)$, and of exponential order p, then
\begin{equation*}
L[(f\ast g)(t)]=L[f(t)].L[g(t)].
\end{equation*}
\end{theorem}
\section{Constructing the propose method}
 In this section, we will investigate solution of fuzzy convolution partial  Volterra integro-differential equation using fuzzy Laplace transform. Consider the most general equation FPIDE
\begin{equation}\label{l1}
\sum_{i=0}^ma_i\frac{\partial^i u}{\partial x^i}+\sum_{i=0}^nb_i\frac{\partial^i u}{\partial t^i}+u+f(x, t)=\int_0^tk(t-s)u(x, s)ds
\end{equation}
With some appropirate fuzzy intial and boundary conditions. Also the functions f(x, t)  is given fuzzy valued function while  k(t, s) is given crisp kernel and $a_i's$, $b_i's$ are constants or functions of x.

\noindent Applying Fuzzy Laplace on both side of  (\ref{l1}) with respect to t we get:
\begin{equation}\label{l2}
\sum_{i=0}^ma_iL[\frac{\partial^i u}{\partial x^i}]+\sum_{i=0}^nb_iL[\frac{\partial^i u}{\partial t^i}]+cL[u]+L[f(x, t)]=L[\int_0^tk(t-s)u(x, s)ds]
\end{equation}
\noindent Using fuzzy  convolution theorem \ref{ct} on integral part we get
\begin{equation}\label{l3}
\sum_{i=0}^ma_iL[\frac{\partial^i u}{\partial x^i}]+\sum_{i=0}^nb_iL[\frac{\partial^i u}{\partial t^i}]+cL[u]+L[f(x, t)]=L[k(t)]L[u(x, t)]
\end{equation}
Now using the definition of fuzzy Laplace stated in section \ref{flt}, equation (\ref{l3}) becomes
\begin{equation}\label{l4}
\sum_{i=0}^ma_i\frac{d^i U(x,p)}{dx^i}+\sum_{i=0}^nb_i U(x,p)\ominus \sum_{j=1}^i p^{j-1}b_i u^{i-j}(x,0)+U(x,p)+L[f(x, t)]=L[k(t)]U(x,p).\end{equation}
The classical form of (\ref{l4}) gives the following two $m^{th}$ order ordinary differential equations as follows:
\begin{equation}\label{l5}
\sum_{i=0}^ma_i\frac{d^i}{dx^i} \underline{U}(x,p,r)+\sum_{i=0}^nb_i \underline{U}(x,p,r)-\sum_{j=1}^i p^{j-1}b_i \underline{u}^{(i-j)}(x,0,r)+\underline{U}(x,p,r)+\underline{F}(x, t,r)=K(t)\underline{U}(x,p,r).\end{equation}
\begin{equation}\label{l6}
\sum_{i=0}^ma_i\frac{d^i}{dx^i} \overline{U}(x,p,r)+\sum_{i=0}^nb_i \overline{U}(x,p,r)-\sum_{j=1}^i p^{j-1}b_i \overline{u}^{(i-j)}(x,0,r)+\overline{U}(x,p,r)+\overline{F}(x, t,r)=K(t)\overline{U}(x,p,r).\end{equation}
Where $K(t)=L[k(t)]$ and $F(x, t)=L[f(x, t)].$ Using the given fuzzy initial and boundary conditions the upper and lower solutions of (\ref{l1}) can be find out from equations (\ref{l5}) and (\ref{l6}) respectively.
\section{Numerical examples}
In this section we will discuss the solution of fuzzy convolution partial  Volterra integro-differential equations using FLT to show the utility of the proposed method in Section 4.
\begin{example}
Let us consider the following fuzzy convolution partial  Volterra integro-differential equation.
\begin{eqnarray}\label{n1}
xu_x=u_{tt}+(x\sin x)(r-1, 1-r)+\int_0^{t}\sin(t-s)u(x,s)ds,
\\ \nonumber \noindent  with \mbox{  } initial \mbox{  } conditions\\ \nonumber \noindent u(x,0,r)=(0,0), \mbox{  } u_t(x,0,r)=((r-1)x, (1-r)x), \\ \nonumber and \mbox{  } Boundary \mbox{  } condition \\ \nonumber u(1, t,r)=((r-1)t, (1-r)t).
\end{eqnarray}
%\subsection{Case $A1$: (i)-differentiability}
Taking  FLT with respect to t  on (\ref{n1}), then we get
\begin{equation}\label{n2}
xL[u_x]=L[u_{tt}]+x(r-1, 1-r)L[\sin t]+L[\sin t]L[u(x,t)].
\end{equation}
Using FLT (\ref{n2}) becomes
\begin{equation}\label{n3}
x\frac{d}{dx}U(x,p)=p^2U(x,p)\ominus p u(x,0)\ominus u_t(x,0)+\frac{x}{p^2+1}(r-1, 1-r)+\frac{1}{p^2+1}U(x,p).
%pL[x(t)]\ominus{x}(0)=(r+1,r-2)L[(1+t)]+L[1]L[x(t)].
\end{equation}
Also applying FLT boundary condition becomes \\
\begin{equation}\label{b1}U(1,p;r)=\frac{(r-1, 1-r)}{p^2}.\end{equation}
The r-cut representation of (\ref{n3}) after using  initials conditions is given by

\begin{equation}\label{n4}
x\frac{d}{dx}\underline{U}(x,p,r)=p^2\underline{U}(x,p,r)-(r-1)x+\frac{x}{p^2+1}(r-1)+\frac{1}{p^2+1}\underline{U}(x,p,r),
\end{equation}
and
\begin{equation}\label{n5}
x\frac{d}{dx}\overline{U}(x,p,r)=p^2\overline{U}(x,p,r)-(1-r)x+\frac{x}{p^2+1}(1-r)+\frac{1}{p^2+1}\overline{U}(x,p,r).
\end{equation}
From (\ref{n4}) we get:
\begin{equation}\label{n6}
\frac{d}{dx}\underline{U}(x,p,r)-\frac{p^4+p^2+1}{x(p^2+1)}\underline{U}(x,p,r)+\frac{(r-1)p^2}{p^2+1}=0,
\end{equation}
Solving (\ref{n6}) we get
\begin{equation}\label{n7}
\underline{U}(x,p,r)=\frac{(r-1)}{p^2}x+Cx^{\frac{p^4+p^2+1}{p^2+1}}.
\end{equation}
On using boundary condition given in (\ref{b1}) we get,  $C=0$ therefore (\ref{n7}) become
\begin{equation}\label{n8}
\underline{U}(x,p,r)=\frac{(r-1)}{p^2}x
\end{equation}
Finally taking inverse Laplace on both side of (\ref{n8})
\begin{equation*}
\underline{u}(x,t,r)=(r-1)xt
\end{equation*}
Similarly on simplifying (\ref{n5}) the following differential equation is obtained:
\begin{equation}\label{n9}
\frac{d}{dx}\overline{U}(x,p,r)-\frac{p^4+p^2+1}{x(p^2+1)}\overline{U}(x,p,r)+\frac{(1-r)p^2}{p^2+1}=0,
\end{equation}
Which gives the final upper  solution of (\ref{n1}) as follows:
\begin{equation*}
\overline{u}(x,t,r)=(1-r)xt
\end{equation*}
\end{example}
\begin{example}
Let us consider the following FPVIDE
\begin{eqnarray}\label{n10}
u_x=u_{tt}+2(1+r,3-r)e^x-2\int_0^{t}(t-s)u(x,s)ds,\\ \nonumber \noindent  with \mbox{  } initial \mbox{  } conditions\\ \nonumber \noindent u_t(x,0,r)=(0,0), \mbox{  } u(x,0,r)=e^x((1+r), (3-r)), \\ \nonumber and \mbox{  } Boundary \mbox{  } condition \\ \nonumber u(0, t, r)=\cos t((r+1), (3-r)).
\end{eqnarray}

Applying FLT on (\ref{n10}), we have
\begin{equation}\label{n11}
L[u_x]=L[u_{tt}]+2(r-1, 1-r)e^xL[1]-2L[t]L[u(x,t)].
\end{equation}
Using definition of FLT (\ref{n11}) becomes
\begin{equation}\label{n12}
\frac{d}{dx}U(x,p)=p^2U(x,p)\ominus p u(x,0)\ominus u_t(x,0)+\frac{2e^x}{p}(r+1, 3-r)-\frac{2}{p^2}U(x,p).
\end{equation}
\noindent After using FLT boundary condition gives:\\
\begin{equation}\label{b2}U(0, p, r)=\frac{(r+1, 3-r)p}{p^2+1}.\end{equation}
The classical form of (\ref{n12}) after using initial conditions, is
\begin{equation}\label{n13}
\frac{d}{dx}\underline{U}(x,p,r)=p^2\underline{U}(x,p,r)-pe^x(r+1)+\frac{2e^x}{p^2+1}(r+1)-\frac{2\underline{U}(x,p,r)}{p^2},
\end{equation}
and
\begin{equation}\label{n14}
\frac{d}{dx}\overline{U}(x,p,r)=p^2\overline{U}(x,p,r)-pe^x(3-r)+\frac{2e^x}{p^2+1}(3-r)-\frac{2}{p^2}\overline{U}(x,p,r),
\end{equation}
Now solving (\ref{n13}) and (\ref{n14}) after using boundary condition (\ref{b2}) we get:
\begin{equation}\label{n15}
\underline{U}(x,p,r)=(r+1)e^x(\frac{p^3-2p}{p^4-p^2-2})
\end{equation}
\begin{equation}\label{n16}
\overline{U}(x,p,r)=(3-r)e^x(\frac{p^3-2p}{p^4-p^2-2})
\end{equation}
Finally taking inverse Laplace we get the lower and upper solutions as follow:
\begin{equation*}
\underline{u}(x,t,r)=(r+1)e^x\cos t
\end{equation*}
\begin{equation*}
\overline{u}(x,t,r)=(3-r)e^x\cos t
\end{equation*}
%\begin{equation}
\end{example}

\begin{example}
Let us consider the following FPVIDE
\begin{eqnarray}\label{n17}
u_{xx}=u_{t}+u+(1+r,3-r)\{-(x^2+1)e^t+2\}+\int_0^{t}e^{(t-s)}u(x,s)ds,\\ \nonumber \noindent  with \mbox{  } given \mbox{  } conditions\\ \nonumber \noindent u(x,0,r)=(r+1,3-r)x^2,  \mbox{  } u_t(x,0,r)=(1+r, 3-r), \\ \nonumber and  \mbox{  }  \\ \nonumber u(0, t, r)=(r+1, 3-r)t, \mbox{  } u_x(0,t,r)=(0, 0).
\end{eqnarray}

Applying FLT on (\ref{n17}), we have
\begin{equation}\label{n18}
L[u_{xx}]=L[u_{t}]+L[u]+(r+1, 3-r)\{-(x^2+1)L[e^t]+2L[1]\}+L[e^t]L[u(x,t)].
\end{equation}
Using definition of FLT (\ref{n18}) becomes
\begin{equation}\label{n19}
\frac{d^2}{dx^2}U(x,p)=pU(x,p)\ominus u(x,0)+U(x,p)+[\frac{-(x^2+1)}{p-1}+\frac{2}{p}](r+1, 3-r)+\frac{U(x,p)}{p-1},
\end{equation}
while the transform boundary conditions are:
\[U(0, p, r)=\frac{(r+1, 3-r)}{p^2}, U_x(0,p,r)=(0, 0)\]\\
The classical form of (\ref{n19}) after using initial conditions, is as under
\begin{equation}\label{n20}
\frac{d^2}{dx^2}\underline{U}(x,p,r)=(\frac{p^2}{p-1})\underline{U}(x,p,r)-(1+r)[x^2+\frac{(x^2+1)}{p-1}]+\frac{2(1+r)}{p},
\end{equation}
and
\begin{equation}\label{n21}
\frac{d^2}{dx^2}\overline{U}(x,p,r)=(\frac{p^2}{p-1})\overline{U}(x,p,r)-(3-r)[x^2+\frac{(x^2+1)}{p-1}]+\frac{2(3-r)}{p},
\end{equation}
Solving (\ref{n20})  we get:
\begin{equation}\label{n22}
\underline{U}(x,p,r)=c_1e^{\sqrt{\frac{p^2}{p-1}}x}+c_2e^{-\sqrt{\frac{p^2}{p-1}}x}+\frac{(1+r)x^2}{p}+\frac{(1+r)}{p^2},
\end{equation}
Using the transform  fuzzy boundary conditions we have $c_1=0$ and $c_2=0$. Hence (\ref{n22}) becomes:
\begin{equation}\label{n23}
\underline{U}(x,p,r)=\frac{(1+r)x^2}{p}+\frac{(1+r)}{p^2}.
\end{equation}
Finally taking inverse Laplace we get
\begin{equation*}
\underline{u}(x,t,r)=(1+r)[x^2+t].
\end{equation*}
Similarly solving (\ref{n21}) in same way we have:
\begin{equation*}
\overline{u}(x,t,r)=(3-r)[x^2+t].
\end{equation*}

\end{example}
\section{Conclusion}
In this paper we investigated the applicability of fuzzy Laplace transform for the solution of FPVIDEs under H-differentiability with crisp kernel. In our knowledge this is the first attempt toward the solution of such equations with fuzzy conditions. We have illustrated the method by solving some examples. In future we will discuss the solution of FPVIDEs under generalized H-differentiability with both crisp and fuzzy kernel.

\bibliography{references1}
\end{document}